\documentstyle[11pt]{article}

\def\@begintheorem#1#2{\par\bgroup{\sc #1\ #2. }\it\ignorespaces}
\def\@opargbegintheorem#1#2#3{\par\bgroup{\sc #1\ #2\ (#3).
}\it\ignorespaces}
\def\@endtheorem{\egroup}
\def\proof{\par{\it Proof}. \ignorespaces}
\def\endproof{{\ \vbox{\hrule\hbox{%
   \vrule height1.3ex\hskip0.8ex\vrule}\hrule
  }}\par}
\newtheorem{theorem}{Theorem}[section]
\newtheorem{lemma}[theorem]{Lemma}

\newcommand{\algSGE}{{2.1}}			%
\newcommand{\algGKOC}{{2.2}}			%
\newcommand{\beq}{\begin{equation}}
\newcommand{\eeq}{\end{equation}}
\newcommand{\beqn}{\[}
\newcommand{\eeqn}{\]}
\newcommand{\bu}{\mbox{\bf u}}
\newcommand{\bv}{\mbox{\bf v}}
\newcommand{\bw}{\mbox{\bf w}}
\newcommand{\bl}{\mbox{\bf l}}
\newcommand{\bs}{\mbox{\bf s}}
\newcommand{\bt}{\mbox{\bf t}}
\newcommand{\bo}{\mbox{\bf 0}}
\newcommand{\ba}{\mbox{\bf a}}
\newcommand{\bb}{\mbox{\bf b}}
\newcommand{\bc}{\mbox{\bf c}}
\newcommand{\bd}{\mbox{\bf d}}
\newcommand{\br}{\mbox{\bf r}}
\newcommand{\be}{\mbox{\bf e}}
\newcommand{\bi}{\mbox{\bf 1}}
\newcommand{\bx}{\mbox{\bf x}}
\newcommand{\bp}{\mbox{\bf p}}
\newcommand{\boldf}{\mbox{\bf f}}
\newcommand{\boldy}{\mbox{\bf y}}
\newcommand{\bomgii}{\mbox{\boldmath$\omega$}_{:2}}
\newcommand{\tomgii}{\tilde{\mbox{\boldmath$\omega$}}_{:2}}
\newcommand{\bgami}{\mbox{\boldmath$\gamma$}_{1:}}

\newcommand{\bmax}{b_{\max}}
\newcommand{\bmin}{b_{\min}}
\newcommand{\bratio}{\frac{\bmax}{\bmin}}
\newcommand{\tlk}{\tilde{\bl}_{:k}}
\newcommand{\tuk}{\tilde{\bu}_{k:}}
\newcommand{\atlk}{\tilde{\bl}_{:k}}
\newcommand{\atuk}{\tilde{\bu}_{k:}}

\newcommand{\lkr}{\bl_{:k}}
\newcommand{\rkr}{\br_{k:n,k}^{(k)}}
\newcommand{\ukr}{\bu_{k:}}
\newcommand{\tlkr}{\tilde{\bl}_{:k}}
\newcommand{\trkr}{\tilde{\br}_{k:n,k}^{(k)}}
\newcommand{\tukr}{\tilde{\bu}_{k:}}

\newcommand{\trkkk}{\tilde{r}_{kk}^{(k)}}

\newcommand{\tbl}{\tilde{\bl}}
\newcommand{\tbr}{\tilde{\br}}
\newcommand{\tbu}{\tilde{\bu}}
\newcommand{\tl}{\tilde{l}}
\newcommand{\tr}{\tilde{r}}
\newcommand{\tu}{\tilde{u}}
\newcommand{\tL}{\tilde{L}}
\newcommand{\tU}{\tilde{U}}
\newcommand{\tR}{\tilde{R}}
\newcommand{\tF}{\tilde{F}}

\newcommand{\Del}[1]{\Delta^{(#1)}}
\newcommand{\del}[1]{\mbox{\boldmath$\delta$}^{(#1)}}
\newcommand{\hDel}[1]{\hat{\Delta}^{(#1)}}
\newcommand{\calD}[1]{{\cal D}^{(#1)}}
\newcommand{\caldi}[1]{\partial_k^{(#1)}}
\newcommand{\phik}{\Phi^{(k)}}
\newcommand{\psik}{\Psi^{(k)}}

\newcommand{\phikk}{\phi_k^{(k)}}
\newcommand{\psikk}{\psi_k^{(k)}}
\newcommand{\phiki}{\Phi^{(k+1)}}
\newcommand{\psiki}{\Psi^{(k+1)}}
\newcommand{\phip}{\Phi^{\prime}}
\newcommand{\psip}{\Psi^{\prime}}
\newcommand{\tphi}{\tilde{\Phi}}
\newcommand{\tpsi}{\tilde{\Psi}}
\newcommand{\tphik}{\tilde{\Phi}^{(k)}}
\newcommand{\tpsik}{\tilde{\Psi}^{(k)}}
\newcommand{\tphiki}{\tilde{\Phi}^{(k+1)}}
\newcommand{\tpsiki}{\tilde{\Psi}^{(k+1)}}
\newcommand{\tphip}{\tilde{\Phi}^{\prime}}
\newcommand{\tpsip}{\tilde{\Psi}^{\prime}}
\newcommand{\tphini}{\tilde{\Phi}^{(n+1)}}
\newcommand{\tpsini}{\tilde{\Psi}^{(n+1)}}

\newcommand{\tphii}{\tilde{\Phi}^{(1)}}
\newcommand{\tpsii}{\tilde{\Psi}^{(1)}}
\newtheorem{corollary}[theorem]{Corollary}
\newcounter{algctr}
\newcommand{\thealgorithm}{\thesection.\arabic{algctr}}
\newcommand{\algorithm}[1]{\addtocounter{algctr}{1}
                          \vspace{2.0mm}\noindent
                          {\bf Algorithm \thealgorithm \ #1}
			  \vspace{-2.5mm}}
\newcommand{\Section}[1]{\section{#1}\setcounter{algctr}{0}}
\newcommand{\Label}[1]{{\rm \label{#1}}}
\newcommand{\for}{{\bf for\ }}
\newcommand{\End}{{\bf end}}
\newcommand{\diag}{{\rm{diag}}}

\begin{document}
\bibliographystyle{plain}
\title{Error Analysis of a Partial Pivoting Method\\
	for Structured Matrices%
\thanks{Copyright \copyright\ 1995, the authors.
	Shorter version to appear in
	{\em Advanced Signal Processing Algorithms},
	Proc.\ SPIE 40th Annual Meeting, San Diego, July 1995.
	\hspace*{\fill} rpb157tr %
        }
}
\author{Douglas R.\ Sweet\\
        Maritime Operations Division\\
        Defence Science and Technology Organisation\\
	Salisbury, SA 5108
\and
	Richard P.\ Brent\\
	Computer Sciences Laboratory\\
	Australian National University\\
	Canberra, ACT 0200}

\date{Report TR-CS-95-03\\
	7 June 1995}

\maketitle				%
\thispagestyle{empty}			%
\vspace{-0.8cm}
\begin{abstract}

\smallskip
Many matrices that arise in the solution of signal processing problems have
a special {\em displacement structure}. For example, adaptive filtering
and direction-of-arrival estimation yield
matrices of Toeplitz type. A recent method of Gohberg, Kailath and Olshevsky (GKO)
allows fast Gaussian elimination with partial pivoting for such structured matrices.
In this paper, a rounding error analysis is performed on the Cauchy and Toeplitz
variants of the GKO method. It is shown the error growth depends on the growth
in certain auxiliary vectors, the {\em generators}, which are computed by the
GKO algorithms. It is also shown that in certain circumstances, the growth in
the generators can be large, and so the error growth is much larger than would
be encountered with normal Gaussian elimination with partial pivoting. A modification
of the algorithm to perform a type of row-column pivoting is proposed which may
ameliorate this problem.
\end{abstract}

\Section{Introduction}
\label{sec:Intro}

Many problems which occur in signal processing, control theory and interpolation
lead to a square or rectangular system with a special structure, for which
the exact or least-squares solution is required. For example,
adaptive filtering requires either the exact solution of a square Toeplitz system
or the least-squares solution of a rectangular Toeplitz system.
A {\em Toeplitz} matrix
is one whose entries along the NW to SE diagonals are constant, i.e. element $t_{ij}$ depends only
on $i-j$. Other types of structured matrices which arise are {\em Hankel} matrices
whose entries along the SW to NE diagonals are constant, {\em Vandermonde} matrices
whose entries have the form $v_{ij} = x_i^{j-1}$, %
and {\em Cauchy} matrices whose entries have the form $c_{ij} = 1/(t_i - s_j)$,
where the $t_i$ and $s_j$ are the elements of vectors {\bf t} and {\bf s}.

Normally, the exact or least-squares solution of a linear system requires $O(n^3)$
operations to solve, where $n$ is the order of the system. However, the structure
of the systems mentioned above 
has been exploited in the past~\cite{TS,rpb108,HVS}
to derive {\em fast} solvers,
i.e. those that require $O(n^2)$ or fewer operations. 
These fast algorithms
are in general numerically unstable for
indefinite systems~\cite{rpb144,Bunch,tpiv}. Recently,
methods have been proposed~\cite{Chan,Gutk,tpiv} which are numerically stable, but which
attempt to retain the $O(n^2)$ complexity. However, all of these algorithms will require
$O(n^3)$ operations in the worst case.
The BBH method~\cite{rpb143}
requires $O(n^2)$ operations in the worst case and
can be shown to be {\em weakly stable}, but not stable in the usual sense of
backward error analysis. Thus there is an interest in fast algorithms which
require $O(n^2)$ operations in the worst case and can be shown to be stable.

Recently, Gohberg, Kailath and Olshevsky~\cite{GKO} have shown how to perform Gaussian
elimination in a fast way with matrices with a special {\em displacement structure}.
Such matrices include Toeplitz, Vandermonde, Hankel and Cauchy matrices, and generalizations
thereof, called {\em Toeplitz type}, etc. They also show how to incorporate partial pivoting
into the Cauchy and Vandermonde solvers. They point out that although pivoting cannot be
incorporated directly into the corresponding Toeplitz or Hankel solvers, the Toeplitz
and Hankel problems can be transformed by simple orthogonal operations into Cauchy problems.
The solution to the original systems can be recovered from those of the transformed
systems by the reverse orthogonal operations. Thus fast Gaussian elimination with partial
pivoting can be carried out on Toeplitz, Vandermonde, Hankel and Cauchy systems.

It might be assumed that such fast solvers should have the same stability properties as
Gaussian elimination with partial pivoting. One of the aims of this paper is analyse the
error behaviour of these algorithms by means of a backward error analysis. 
It is shown
that error propagation depends on the magnitude of both the triangular factors
$L$ and $U$ (as in Gaussian elimination) and the {\em generators}, auxiliary vectors which
are computed during the course of the algorithm.

It is shown that in some cases the generators can suffer a large growth and cause a
corresponding growth in the backward and forward error.
A modification is proposed
which may prevent this growth,
and so restore the stability of the algorithm in these cases.
However, we can not prove that the modification is always successful.

The paper is structured as follows. In \S\ref{sec:GKO}, the
Gohberg-Kailath-Olshevsky (GKO) algorithm for Cauchy and
Toeplitz matrices is briefly described. The error analyses of the Cauchy and Toeplitz
variants of the GKO algorithm are carried out in \S\ref{sec:GKOerror}
and \S\ref{sec:GKOTerror} respectively.
In \S\ref{sec:Discussion}, examples for both variants are given where
a large growth occurs in the generators and hence in
the errors in the solutions.
The modified version of the GKO algorithm is proposed in \S\ref{sec:GKOmod},
and numerical tests of this are carried out there. Some conclusions are
drawn and suggestions for future work are given in \S\ref{sec:Conclusion}.

{\em Notation.}
\ The following notation is used. $\epsilon$ is the machine epsilon,
and $n$ is the order of the matrix to be factorized. Scalars of the form
$c_i$ and $k_i$ are small constants. $\be_j$ denotes the $j$th
column of the identity matrix. Elementwise matrix multiplication is
denoted by the centred circle $\circ$. For a matrix $A$, $|A|$ is the
matrix of moduli of the $\{a_{ij}\}$, $A^I$ denotes elementwise inversion,
and $A^{\prime}$ denotes augmentation of $A$ to order $n$ by adding zero rows
and zero columns respectively above and to the left of $A$.
Other submatrices are indicated in MATLAB
style, i.e. for a matrix $A$, $A_{p:q,r:s}$ selects rows $p$ to $q$ of
columns $r$ to $s$, and a colon without an index range selects all of the
rows or columns.

\pagebreak[4]
\Section{The Gohberg-Kailath-Olshevsky (GKO) \mbox{Algorithm}}
\label{sec:GKO}

In this section, we first define the displacement operator, displacement
equation and displacement rank for structured matrices; we then give the
general Gaussian elimination algorithm for structured matrices, followed
by the variants for Cauchy and Toeplitz matrices.

\subsection{Displacement structure}
Gohberg {\em et al}\,~\cite {GKO} show that structured matrices satisfy a
{\em Sylvester equation} which has the form
   \beq \nabla_{\{A_f,A_b\}}(R) = A_f R - R A_b = \Phi\Psi\;, \label{sylv} \eeq 
where $A_f$ and $A_b$ have some simple structure (usually banded, with 3 or fewer
full diagonals), $\Phi$ and $\Psi$ are $n \times \alpha$ and $\alpha \times n$ respectively,
and $\alpha$ is some small integer (usually 4 or less). The pair of matrices
$\Phi, \Psi$ is called
the $\{A_f, A_b\}$-{\em generator} of $R$, and $\alpha$ is called the 
$\{A_f, A_b\}$-{\em displacement rank} of $R$.

Particular choices of $A_f$ and $A_b$ lead to definitions of basic classes of 
matrices. Thus, for a Cauchy matrix 
$$C(\bt,\bs) = \left[ \frac{1}{t_i - s_j} \right]_{ij}\;,$$
we have
   \beq A_f = D_t = \mbox{diag}(t_1, t_2, \ldots, t_n), \quad 
     A_b = D_s = \mbox{diag}(s_1, s_2, \ldots, s_n)  \label{Af.Cauchy} \eeq
and
   \beqn \Phi^T = \Psi = [1,1,\ldots,1]\;. \eeqn
More general matrices, where $A_f$ and $A_b$ are as in (\ref{Af.Cauchy}) but 
$\Phi$ and
$\Psi$ are general rank-$\alpha$ matrices, are called {\em Cauchy-type}.

Similarly, for a Toeplitz matrix $T = [t_{ij}] = [a_{i-j}]$
\beq  A_f = Z_1 = \left[ \begin{array}{ccccc}
                         0 &   0    & \cdots & 0 & 1      \\
                         1 &   0    &        &   & 0      \\
                         0 &   1    &        &   & \vdots \\
                    \vdots &        & \ddots &   & \vdots \\
                         0 & \cdots &    0   & 1 & 0
                         \end{array}  \right], \hspace{5mm}
      A_b = Z_{-1} =  \left[ \begin{array}{ccccc}
                         0 &   0    & \cdots & 0 &  -1    \\
                         1 &   0    &        &   &   0    \\
                         0 &   1    &        &   & \vdots \\
                    \vdots &        & \ddots &   & \vdots \\
                         0 & \cdots &    0   & 1 &   0
                         \end{array}  \right],   \label{Af.Toep} \eeq
\beq \Phi = \left[ \begin{array}{ccccc}
                   1  &    0    &   \cdots   &    \cdots      &      0        \\
                  a_0 & a_{1-n}+a_1 & \cdots & a_{-2}+a_{n-2} & a_{-1}+a_{n-1}
                \end{array} \right]^T            \label{Phi.Toep}\eeq
and
\beq \Psi =\left[ \begin{array}{ccccc}
                a_{n-1}-a_{-1} & a_{n-2}-a_{-2} & \cdots & a_1-a_{1-n} & a_0 \\
                       0       &     \cdots     & \cdots &   0    &  1 %
                \end{array} \right]\;.            \label{Psi.Toep}\eeq

\subsection{Gaussian elimination for structured matrices}
Let the input matrix, $R_1$, have the partitioning 
$R_1 = \left[ \begin{array}{cc} d_1 & \bw_1^T \\ \boldy_1 & \dot{R}_1
       \end{array} \right]$.
The first step of normal Gaussian elimination is to premultiply $R_1$ by
      $\left[ \begin{array}{cc} 1 & \bo^T \\ -\boldy_1/d_1 & I
       \end{array} \right]$,
which reduces $R_1$ to 
      $\left[ \begin{array}{cc} d_1 & \bw_1^T \\ \bo & R_2
       \end{array} \right]$,
where $R_2 = \dot{R}_1 - \boldy_1\bw_1^T/d_1$ is the {\em Schur
complement} of $d_1$ in $R_1$. At this stage, $R_1$ has the factorization
\beqn R_1 = \left[ \begin{array}{cc} 1 & \bo^T \\ \boldy_1/d_1 & I
       \end{array} \right]
       \left[ \begin{array}{cc} d_1 & \bw_1^T \\ \bo & R_2
       \end{array} \right]\;. \eeqn
One then proceeds recursively with the Schur complement
$R_2 = \left[ \begin{array}{cc} d_2 & \bw_2^T \\ \boldy_2 & \dot R_2
       \end{array} \right]$, eventually yielding a factorization $R_1 = LU$,
where column $k$ of $L$ is $[\bo^T \quad 1 \quad \boldy_k^T]^T$, and row
$k$ of $U$ is $[\bo^T \quad 1 \quad \bw_k^T]$.

The genesis of {\em structured} Gaussian elimination is the fact that
the displacement structure is preserved under Schur complementation, and
that the generators for the Schur complement $R_{k+1}$ can be computed from
the generators of $R_k$ in $O(n)$ %
operations. This is expressed constructively
in the following theorem, which is proved in~\cite{GKO}.
\begin{theorem}
Let matrix $R_1 = \left[ \begin{array}{cc} d_1 & \bw_1^T \\ \boldy_1 & \dot R_1
                         \end{array} \right]$
satisfy the Sylvester equation
   \beq \nabla_{\{A_{f,1},A_{b,1}\}}(R_1) = A_{f,1} R_1 - R_1 A_{b,1} = 
   \Phi^{(1)}\Psi^{(1)}\;,  \label{Sylv1} \eeq
where $\Phi^{(1)} = [\varphi_1^{(1)T} \quad \varphi_2^{(1)T} \quad \cdots
                                \quad \varphi_n^{(1)T}]^T $,
$\Psi^{(1)} = [\psi_1^{(1)} \quad \psi_2^{(1)} \quad \cdots
                            \quad \psi_n^{(1)}]$,
$\varphi_i^{(1)} \in \mbox{\bf C}^{1\times\alpha}$ and
$\psi_i^{(1)} \in \mbox{\bf C}^{1\times\alpha}$, $(i = 1, 2, \ldots, n)$.
Then $R_2$, the Schur complement of $d_1$ in $R_1$, satisfies the Sylvester
equation
   \beqn \nabla_{\{A_{f,2},A_{b,2}\}}(R_2) = A_{f,2} R_2 - R_2 A_{b,2} = 
   \Phi^{(2)}\Psi^{(2)}\;,  \eeqn
where $A_{f,2}$ and $A_{b,2}$ are respectively $A_{f,1}$ and $A_{b,1}$ 
with their first rows and first columns deleted, and 
where $\Phi^{(2)} = [0,\; \varphi_2^{(2)T},\; \varphi_3^{(2)T},\; \cdots
                                \;, \varphi_n^{(2)T}]^T $
and $\Psi^{(2)} = [0,\; \psi_2^{(2)},\; \psi_3^{(2)},\; \cdots
                            \;, \psi_n^{(2)}]$ are given by 
\beq \Phi_{2:n,:}^{(2)} = 
     \Phi_{2:n,:}^{(1)} - \boldy_1\varphi_1^{(1)}/d_1 \label{phi.ud}\;,\eeq
\beq \Psi_{:,2:n}^{(2)} = 
     \Psi_{:,2:n}^{(1)} - \psi_1^{(1)}\bw_1^T/d_1 \label{psi.ud}\;. \eeq
\end{theorem}

Equations (\ref{phi.ud}) and (\ref{psi.ud}) form the basis of the following
general structured Gaussian elimination algorithm.

\algorithm{(Structured Gaussian elimination)}
\begin{enumerate}
   \item Recover from the generator $\Phi^{(1)}$, $\Psi^{(1)}$ the first row
         and column of
         $R_1 = \left[ \begin{array}{cc} d_1 & \bw_1^T \\ \boldy_1 & R_{22}^{(1)}
                         \end{array} \right]\;.$
   \item $[1 \quad \boldy_1^T/d_1]^T$ and $[d_1 \quad \bw_1^T]$ are respectively
         the first column and row of $L_1$ and $U_1$ in the $LU$ factorization
         of $R_1$.
   \item Compute by equations (\ref{phi.ud}) and (\ref{psi.ud}),
         the generator $\Phi^{(2)}, \Psi^{(2)}$ for
         the Schur complement $R_2$.
   \item Proceed recursively with $\Phi^{(2)}$ and $\Psi^{(2)}$. Each major
         step yields $[1 \quad \boldy_k^T/d_k]^T$ and $[d_k \quad \bw_k^T]$,
         which are respectively
         the first column and row of $L_k$ and $U_k$ in the $LU$ factorization
         of $R_k$. Column $k$ of $L$ and row $k$ of $U$ are respectively
         $[\bo_{k-1}^T \quad 1 \quad \boldy_k^T/d_k]^T$ and
         $[\bo_{k-1}^T \quad d_k \quad \bw_k^T]$.
\end{enumerate}

{\bf Pivoting.} Gaussian elimination without pivoting is unstable in general.
One normally uses partial pivoting (swapping rows to bring the largest
element in the first column to the pivot position) or complete pivoting (swapping
rows and columns to bring the largest element in the whole matrix to the
pivot position) to improve the accuracy. Row and/or column interchanges
can destroy the structure of certain matrices, such as Toeplitz matrices.
However, if $A_{f,1}$ in (\ref{Sylv1}) is diagonal (which is the case for
Cauchy and Vandermonde type matrices), then the structure is preserved
under row permutations.

Partial pivoting can also be incorporated into structured Gaussian elimination.
Suppose we wish to swap rows 1 and $q$ of $R_1$. Let $P_1$ be
the matrix which applies this permutation. Then it is easy to see that
$P_1 R_1$ satisfies (\ref{Sylv1}) with the $(1,1)$ and $(q,q)$ entries of
$A_{f,1}$ swapped, and with swapped row vectors $\varphi_1^{(1)}$ and
$\varphi_q^{(1)}$. Thus, pivoting can be incorporated into
Algorithm~\thealgorithm\, by adding the following steps:

\pagebreak[3]
\begin{itemize}
   \item[0.5] {\em Initialization step.} Set permutation matrix $P = I$.
   \item[2.5] {\em After step 2 of Algorithm~\thealgorithm.}
          Let $(\boldy_1)_q$ be the largest
          entry by magnitude in $\boldy_1$. Swap rows $1$ and $q$ of $P$
          and $\Phi^{(1)}$, and the first and $q$-th diagonal entries
          in $A_{f,1}$. Recover the first row of $P_1 R_1$ from
          $\Psi^{(1)}$ and the swapped $\Phi^{(1)}$.
\end{itemize}
{\em Note.} The computation of the first row of the original $R^{(1)}$ in
step 1 of Algorithm~\thealgorithm\, may be omitted~-- we only require the
first row of the {\em swapped} $R^{(1)}$.

It may be seen that the pivoted algorithm computes upper and lower triangular
matrices $L$ and $U$ which satisfy
\beqn R^{(1)} = P^T L U\;.\eeqn
Note that for Cauchy-type matrices, where both $A_{f,1}$ and $A_{b,1}$ are
diagonal, both row and column pivoting may be performed. However, complete
pivoting requires the computation of all the entries in the matrix, which
would require $O(n^2)$ operations at each step and $O(n^3)$ operations in all.
It will be seen in \S\ref{sec:GKOmod} that a restricted version of row-column
pivoting can be used to improve the performance of the GKO algorithm.

\subsection{The Cauchy variant of the GKO algorithm (GKO-Cauchy)}
Recall that a Cauchy-type matrix satisfies the Sylvester equation (\ref{Sylv1})
with 
   $$ A_{f,1} = D_t = \mbox{diag}(t_1, t_2, \ldots, t_n) \quad \mbox{and} \quad 
     A_{b,1} = D_s = \mbox{diag}(s_1, s_2, \ldots, s_n)\;. $$
It can be easily verified that if $t_i \neq s_j$, then the $(i,j)$ entry of
$R^{(1)} = R$ is given by
\beqn r_{ij} = \frac{\varphi_i \psi_j}{t_i - s_j}\;.\eeqn
There may be some cases where $t_i = s_j$ and $\varphi_i \psi_j = 0$ for some
$(i,j)$, and $r_{ij}$ cannot be recovered from its generator. We do not
consider these cases in this paper.

In general, at major step $k$, the reduced matrix $R^{(k)}$ has the form
$$ R^{(k)} = \left[ \begin{array}{cccc} 
                      d_1  &        &         & \bw_1^T     \\
                      \bo  & \ddots &         & \vdots      \\
                    \vdots &        & d_{k-1} & \bw_{k-1}^T \\
                      \bo  & \cdots &   \bo   & R_k
                    \end{array} \right]\;. $$
The entries of the $k$-th Schur component $R_k$, may be computed by
\begin{eqnarray}
     r_{ij}^{(k)} & = & \frac{\varphi_i^{(k)} \psi_j^{(k)}} {t_i - s_j}\;,
                        \quad k \leq i,j \leq n \label{rgen} \\
                  & = & (R_k)_{i-k+1,j-k+1} \nonumber
\end{eqnarray}
Equation~(\ref{rgen}) can be used in Algorithm~\thealgorithm\, with pivoting
to yield the Cauchy version of the GKO algorithm.

\pagebreak[3]
\algorithm{(GKO-Cauchy)}
\begin{tabbing}
12345\=123\=123\=123456789012345678901234567890\=\kill \\
{\em Input.}  \>\> Cauchy-type matrix $R_1$, specified by $\bt$, $\bs$,
                   $\Phi^{(1)}$ and $\Psi^{(1)}$. \\
{\em Output.} \>\> Factorization $R_1 = P^TLU$, where $P$ is a permutation,
                   and $L$ and \\
              \>\> $U$ are lower and upper-triangular respectively. \\
\\
\> \% Initialization \\
\> $L \leftarrow 0 ; \quad U \leftarrow 0 ; \quad P \leftarrow I$ \\
\> \for $k \leftarrow 1:n$  \>\>\> \% $k$: Iteration number  \\
\> \> \for $j \leftarrow k:n$ \>\> \% recover col.1 of $R_k$ \\
\> \> \> $r_{jk}^{(k)} \leftarrow \frac{\varphi_j^{(k)}\psi_k^{(k)}}{t_j - s_k}$ \\
\> \> \End \\
\> \> \% Carry out row interchanges     \\
\> \> Find $k \leq q \leq n$ such that
      $|r_{qk}^{(k)}| = \max_{k\leq j\leq n} |r_{jk}^{(k)}|$ \\
\> \> $t_k \leftrightarrow t_q ; \quad
      \varphi_k^{(k)} \leftrightarrow \varphi_q^{(k)} ; \quad
      r_{kk}^{(k)} \leftrightarrow r_{qk}^{(k)}$         \\
\> \> {\bf swap} $k$-th and $q$-th rows of $L$             \\
\> \> {\bf swap} $k$-th and $q$-th rows of $P$             \\
\\
\> \> \for $j \leftarrow k+1:n$ \>\> \% Recover row 1 of swapped $R_k$ \\
\> \> \> $r_{kj}^{(k)} \leftarrow \frac{\varphi_k^{(k)}\psi_j^{(k)}}{t_k - s_j}$ \\
\> \> \End \\
\> \> $u_{kk} \leftarrow r_{kk}^{(k)}$ \\
\> \> \% Compute row and col $k$ of $L$ and $U$, and update $\Phi$ and $\Psi$
         using (\ref{phi.ud}) and (\ref{psi.ud}) \\
\> \> \for $j \leftarrow k+1:n$ \\
\> \> \> $l_{jk} \leftarrow r_{jk}^{(k)}/r_{kk}^{(k)}$                      \\
\> \> \> $u_{kj} \leftarrow r_{kj}^{(k)}$                                   \\
\> \> \> $\psi_j^{(k+1)} \leftarrow \psi_j^{(k)} - \psi_k^{(k)}u_{kj}/u_{kk}$ \\
\> \> \> $\varphi_j^{(k+1)} \leftarrow \varphi_j^{(k)} - \varphi_k^{(k)}l_{jk}$ \\
\> \> \End   \\
\> \End      \\
\end{tabbing}

\subsection{The Toeplitz variant of the GKO algorithm (GKO-Toeplitz)}
Recall that a Toeplitz matrix satisfies the Sylvester equation (\ref{sylv}),
with $A_f$, $A_b$, $\Phi$ and $\Psi$ being given by 
equations~(\ref{Af.Toep}) to
(\ref{Psi.Toep}), and a Toeplitz-type matrix is one with $A_f$ and $A_b$
given by (\ref{Af.Toep}), and with general low-rank $\Phi$ and $\Psi$. The first
row and column of the Toeplitz-type matrix can be simply generated from
$\Phi$ and $\Psi$, and this generating formula can be used in
Algorithm~\algSGE\ to 
yield a structured Gaussian elimination algorithm for Toeplitz-type matrices.

Because neither $A_f$ nor $A_b$ is diagonal, pivoting cannot be introduced directly into this
structured algorithm~-- pivoting will destroy the Toeplitz-type property.
However, the Toeplitz-type matrix can be easily converted, by fast orthogonal
transformations, into a Cauchy-type matrix which can be factorized as in
Algorithm~\algSGE.
The inverse orthogonal transforms yield the factorization
of the original matrix. The following result of~\cite{GKO} shows
how this conversion may be done.

\pagebreak[4]
\begin{theorem}
Let $T$ be a Toeplitz-type matrix, satisfying
$$ \nabla_{\{Z_1,Z_{-1}\}}(T) = \Omega\Gamma\;, $$
$$ \Omega = [\omega_1^T \quad \omega_2^T \quad \cdots \omega_n^T]^T, \quad
   \Gamma = [\gamma_1 \quad \gamma_2 \quad \cdots \quad \gamma_n], $$
where the $\{\omega_i\}$ and the $\{\gamma_i\}$ are $1\times\alpha$ and $\alpha\times1$
respectively. \\
Then 
\beq
R = FTD^{-1}F^\ast \label{T.C} 
\eeq
is a Cauchy-type matrix, satisfying
$$ \nabla_{\{D_F,D_{F\_}\}} = \Phi\Psi\;, $$
where $F = \frac{1}{\sqrt n} [ e^{2\pi i(k-1)(j-1)/n} ]_{1\leq k,j\leq n}$ 
is the Discrete Fourier Transform matrix,
\beq
D_F = \diag(1,e^{2\pi i/n},\ldots,e^{2\pi i(n-1)/n})\;, \quad
D_{F\_} =  \diag(e^{\pi i/n},e^{3\pi i/n},\ldots,e^{\pi i(2n-1)/n}) \;,
                \label{eq:Df}
\eeq
$$
D = \diag(1,e^{\pi i/n},\ldots,e^{\pi i(n-1)/n})   
$$
and
\beq 
\Phi = F\Omega\;, \quad\quad \Psi^\ast = FD\Gamma^\ast\;. \label{T->C} 
\eeq
\label{Tgen.Cgen}
\end{theorem}

Theorem~\ref{Tgen.Cgen} allows the generators of $T$ to be converted to the
generators of $R$ in $O(2\alpha n\log n)$ operations via FFTs. $R$ can then
be factorized as $R = P^TLU$ using Algorithm~\algGKOC.
Using (\ref{T.C}), we then obtain
\beq
T = F^\ast P^T L U F D  \;. \label{T_fac}
\eeq
{From} this factorization, 
a linear system in $T$ can be solved in $n^2 + 2n\log n$
operations, so the whole procedure of conversion of Cauchy form, factorization
and solution requires $O(n^2)$ operations.

\Section{Error Analysis of the GKO-Cauchy Algorithm}
\label{sec:GKOerror}

In this section, a backward error analysis will be carried out, which yields
a bound for the perturbation matrix $E$, defined by
\beq  \tL\tU = R + E\;, \label{E_def} \eeq
where $R$ is the matrix to be factorized, and $\tL$ and $\tU$ are the
{\em computed} factors. In the analysis, we first derive some preliminary
results which apply to {\em any} algorithm for structured Gaussian elimination
(SGE), and indicate a general methodology for error analysis of SGE algorithms.
We then carry out the analysis for Cauchy-type matrices in general and for the
Cauchy-type matrix derived from a Toeplitz matrix by equation~(\ref{T.C}). 

\subsection{Preliminary results}
The following two lemmas may be used for the error analysis of SGE algorithms
in general, and the GKO-Cauchy algorithm in particular. The first lemma
shows that if $G$ is the perturbation in the Sylvester equation caused
by replacing $R$ by $\tL\tU$, then the displacement of $E$ is $G$.
\begin{lemma}
\label{lem1}
Let $R$ be a general structured matrix that satisfies {\rm (\ref{sylv})}, let
$A_f$, $A_b$, $\Phi$ and $\Psi$ be as defined above, and let $\tL$, $\tU$
and $E$ be as in {\rm (\ref{E_def})}. Suppose $\tL$ and $\tU$ satisfy
\beq
A_f\tL\tU - \tL\tU A_b = \Psi\Phi + G \;;   \label{lem1.1}
\eeq
then $E$ satisfies
\beq 
\nabla_{\{A_f,A_b\}}(E) \equiv A_f E - E A_b = G\;. \label{lem1.2} 
\eeq
\end{lemma}
\proof From (\ref{E_def}) and (\ref{lem1.1}),
\beqn
A_f(R+E) - (R+E)A_b = \Phi\Psi + G \;.
\eeqn
Expanding the above, and using (\ref{sylv}) we obtain (\ref{lem1.2}).
\endproof

\begin{corollary}
\label{coroll1}
If $R$ is a Cauchy-type matrix with $A_f = D_t$ and $A_b = D_s$, then $E$
satisfies
\beq
D_t E - ED_s = G   \label{coroll1.1}
\eeq
and
\beq
e_{ij} = \frac{g_{ij}}{t_i-s_j} \;, \quad i,j = 1,\ldots,n  \label{coroll1.2}
\eeq
\end{corollary}
\proof (\ref{coroll1.1}) follows directly from (\ref{lem1.2}),
and (\ref{coroll1.2}) follows by evaluating each component of 
(\ref{coroll1.1}). 
\endproof

If $R$ is a Toeplitz-type matrix, $A_f = Z_1$, $A_b = Z_{-1}$, and $E$ satisfies
$Z_1 E - EZ_{-1} = G$. Because $Z_1$ and $Z_{-1}$ are not diagonal, the
recovery formula for $E$ is a little more involved, and will be derived
in the next section (Lemma~\ref{lem4.2}).

The second lemma of this section shows that $G$ is the sum of the local perturbation
matrices incurred in each step of the relevant structured Gaussian elimination
(SGE) algorithm.
\begin{lemma}
\label{lem2}
Let $\nabla_{\{A_f,A_b\}}$ be the displacement operator as defined
in {\rm (\ref{sylv})}; let $\tL$, $\tU$ and $G$ be as defined above;
let the $\{\tphik, \tpsik\}_{k=1,2,\ldots}$ be the computed generators
of the $\{R_k^{\prime}\}_{k=1,2,\ldots}$, the reduced matrices at step $k$ of
SGE, and define $\tphini = \tpsini = \bo$. Then
\beq
G = \sum_{k=1}^n H_k  \;,     \label{lem2.1}
\eeq
where $H_k$, the local perturbation in each step of SGE, is defined by
\beq
\nabla_{\{A_f,A_b\}}(\tlk\tuk) = 
             \tphik\tpsik - \tphiki\tpsiki + H_k \;,
                                       \quad k = 1,\ldots,n \;. \label{lem2.2}
\eeq
\end{lemma}
\proof Writing (\ref{lem2.2}) explicitly, we get
\beq
A_f\tlk\tuk - \tlk\tuk A_b = \tphik\tpsik - \tphiki\tpsiki + H_k \;,
                                       \quad k = 1,\ldots,n \;. \label{lem2.3}
\eeq
Summing the members of (\ref{lem2.3}), we obtain
\beq
A_f\sum_{i=1}^n\atlk\atuk - \sum_{i=1}^n\atlk\atuk A_b = 
 \tphii\tpsii - \tphini\tpsini + \sum_{i=1}^nH_k \;. \label{lem2.4}
\eeq
Now
$\sum_{i=1}^n\tlk\tuk = \tL\tU$, $\tphii = \Phi$, $\tpsii = \Psi$
and $\tphini \equiv \tpsini \equiv \bo$. Substituting these identities into
(\ref{lem2.4}) and comparing the resulting relation with (\ref{lem1.1}),
we obtain (\ref{lem2.1})
\endproof

\subsection{Methodology of error analysis for SGE algorithms}
Lemmas \ref{lem1} and \ref{lem2} may be used in a general methodology for
the error analysis of SGE algorithms similar to Algorithm~\algSGE.

In the following methodology and the subsequent analysis of the GKO algorithm,
we now let $\phik$ and $\psik$ be the {\em computed} values of these quantities,
$\ukr$, $\rkr$, $\lkr$, $\phiki$ and $\psiki$ be the values of these quantities
computed in exact arithmetic from $\phik$ and $\psik$ using steps 1 to 3 of
Algorithm~\algSGE, and $\tukr$, $\trkr$, $\tlkr$, $\tphiki$ and $\tpsiki$ be the
actual computed values of $\ukr$, $\rkr$, $\lkr$, $\phiki$ and $\psiki$ respectively.
The methodology is as follows:

\pagebreak[3]
\begin{enumerate}
\item Using a standard rounding error analysis, derive expressions of the form
      \begin{eqnarray}
      \tukr & = & \ukr + \delta\tukr   \label{err1}  \\
      \trkr & = & \rkr + \delta\trkr   \label{err2}  \\
      \tlkr & = & \lkr + \delta\tlkr   \label{err3}  \\
      \tphiki & = & \Phi^{(k)} - \tlkr\phikk + \delta\tphiki \label{err4}  \\
      \tpsiki & = & \Psi^{(k)} - \psikk\tukr/\trkkk + \delta\tpsiki       
                                            \label{err5} 
      \end{eqnarray}
      where $\delta\tukr$, etc. are error terms.
\item Evaluate $\phik\psik - \tphiki\tpsiki$ using (\ref{err1}) to (\ref{err5}).
      This can be expressed in the form
      \beq
      \phik\psik - \tphiki\tpsiki = A_f\tlkr\tukr - \tlkr\tukr A_b + F_k \;,   
                                                     \label{F_k}
      \eeq
      where $F_k$ is an error term. By (\ref{lem2.3}),
      $$H_k = -F_k \;.$$ 
\item After some manipulation, express $F_k$ as a sum of terms of
      the form
      \[	
      S(A_f,A_b) \circ T(V^{(k)}) \circ \tlkr\tukr \circ \hat\Delta
	\;\;{\rm or}\;\;
      S(A_f,A_b) \circ T(V^{(k+1)}) \circ L_{:,k+1:n}U_{k+1:n,:}
      \circ \hat\Delta\;.
      \]
      Here, the $S(A_f,A_b)$ are matrices formed from $A_f$ and $A_b$, $\Delta$
      is a matrix whose elements are bounded in magnitude by $\epsilon$, and
      $V^{(k)}$ is defined by 
      \beq
      |\phik||\psik| \equiv V^{(k)} \circ \phik\psik \;.  \label{def_V}
      \eeq
\item Apply (\ref{lem2.1}) to derive an expression for $G$.
\item Lemma~\ref{lem1} shows that $G$ satisfies
      \beq \nabla_{\{A_f,A_b\}} E = G \;.  \label{dispE}  \eeq
      Using the appropriate algorithm to recover a structured matrix from its
      generators, derive an expression for $E$ from the expression for $G$.
      Note that in general, $G$ will be of full rank. However, (\ref{dispE})
      will still be satisfied by $E$ and $G$.
\item Derive bounds for $\|E\|$ using some norm.
\end{enumerate}

\subsection{Error analysis of GKO for Cauchy-type matrices}
In this subsection, we use the above methodology to derive the first of our main
results --- a bound for $\|E\|$ when a Cauchy matrix $R$ is factorized by the
GKO algorithm. The results are encapsulated in three theorems, which yield
expressions for the $\{H_k\}$, an elementwise bound for $G$, and a bound
for $\|E\|$ respectively. We then discuss the size of the bound for $\|E\|$.
\begin{theorem}
\label{Th3.1}
Let $R$ be a Cauchy matrix to be factorized by the GKO algorithm and
let $F_k$, $H_k$, $V^{(k)}$, $\tlkr$, $\tukr$ be as defined above. Then
\begin{eqnarray}
F_k & = & c_1\hDel{1} \circ D_{vc}^{(k)} D_p \tlkr\tukr
 + c_2\tlkr\tukr D_q D_{vr}^{(k)} \circ \hDel{2}
 + c_3 (r_{kk}^{(k)})^{-1} v_{kk}^{(k)} \hDel{3} \circ \tlkr\tukr + \nonumber \\
& & c_4\hDel{4} \circ B^I \circ V^{(k+1)} \circ \tL_{:,k+1:n}\tU_{k+1:n,:}
                  \nonumber %
\end{eqnarray}
where $c_1$ to $c_4$ are small constants, 
$D_{vc}^{(k)} = \diag(v_{:k}^{(k)})$,
$D_{vr}^{(k)} = \diag(v_{k:}^{(k)})$, $D_p = \diag\{t_i-s_k\}_i$,
$D_q = \diag\{t_k-s_j\}_j$, 
$B = [1/(t_i-s_j)]$ is the ordinary Cauchy
matrix with displacement operator $\nabla_{\{D_s,D_t\}}$,
the $\hDel{\cdot}$ are matrices whose elements are less than $\epsilon$
in magnitude, and $H_k = -F_k$.
\end{theorem}
\proof In the following, we simplify our notation and drop the superscript $(k)$;
where the superscript is $(k+1)$ we indicate this by a prime $(^\prime)$; 
and we drop the subscripts $:k$, $k:$ and $k:n,k$.
In the following, we do not give all the steps in the derivation of the
various expressions, as these are straightforward but very tedious.
However, we indicate how key intermediate expressions are derived.

We use the normal properties of
floating point operations performed with at least one guard digit, 
viz.\ $fl(a) = a(1+\delta_1)$ and 
$fl(a\star b) = (a\star b)(1+\delta_2)$, where $fl(a)$ denotes rounding,
$fl(a\star b)$ is the computed result of
any of the four basic floating-point operations, and
$|\delta_1|, |\delta_2| < \epsilon$.

Following step 1 of the above methodology, we evaluate expressions for the
computed values of $\tbr$, $\tbl$ and $\tbu$ (subscripts and superscripts
dropped), yielding after a few steps
\begin{eqnarray}
\tbu & = & \bu + 2\tbu\Del{1} + \phi_k\calD{1}\Psi D_q^{-1} \;, \label{tbu} \\
\tbr & = & \br + 2\Del{2}\tbr + D_p^{-1}\calD{2}\Phi\psi_k  \;,\label{tbr} \\
\tbl & = & \bl + 5\Del{3}\tbl + 
	\tr_{kk}^{-1} D_p^{-1}\calD{2}\Phi\psi_k  \nonumber %
              - b_{kk}\tr_{kk}^{-1}\caldi{2}\phi_k\psi_k\tbl \;.
\end{eqnarray}
Here and below the $\Del{\cdot}$ denote diagonal matrices with
elements of magnitude less than $\epsilon$; the $\calD{\cdot}$ are elementwise
operators which multiply each element of their matrix operands by a factor
less than $\epsilon$, and the $\caldi{\cdot}$ are
similar elementwise vector operators.

Similarly, it can be shown that the computed values of $\phip$ and $\psip$
satisfy
\begin{eqnarray}
\tphip & = & \Phi - \tbl\phi_k + \calD{3}\phip + 
	\calD{4}(\tbl\phi_k) \;, \nonumber \\ %
\tpsip & = & \Psi - \psi_k\tbu/\tr_{kk} + \calD{5}\psip +
        2\calD{6}(\psi_k\tbu)/\tr_{kk} \;.  \nonumber %
\end{eqnarray}

Carrying out step 2 of the above methodology, we obtain
\begin{eqnarray}
\Phi\Psi - \tphip\tpsip & = &
\Phi\psi_k\tbu/\tr_{kk} + \tbl\phi_k\Psi - \tbl\phi_k\psi_k\tbu/\tr_{kk}
- 2\Phi\calD{6}(\psi_k\tbu)/\tr_{kk} - \phip\calD{5}\psip - 
    \nonumber \\
&  &
\calD{4}(\tbl\phi_k)\Psi - \calD{3}\phip\psip +
\calD{4}(\tbl\phi_k)\psi_k\tbu/\tr_{kk}
+ 2\tbl\phi_k\calD{6}(\psi_k\tbu)/\tr_{kk}  \;.     \label{11_terms}
\end{eqnarray}
Let $T_3$ denote the first three terms in (\ref{11_terms}).
{From} Algorithm~\algGKOC,
we have $\Phi\psi_k = D_p\br$ and $\phi_k\Psi = \tbu D_q$.
Using these relations in $T_3$, and expressing $\br$ in terms of ($\tbr$ - error terms)
using (\ref{tbr}) and $\bu$ in terms of ($\tbu$ - error terms) using (\ref{tbu}),
we can show that
\begin{eqnarray}
T_3 & = & D_t\tbl\tbu - \tbl\tbu D_s - 3D_p\Del{4}\tbl\tbu - 2\tbl\tbu\Del{5}D_q
+ 2r_{kk}^{-1}\delta\tbl\tbu -
\nonumber \\
& &
\calD{2}\Phi\psi_k\tbu/r_{kk} - \tbl\phi_k\calD{1}\Psi +
\tr_{kk}^{-1}\caldi{2}\phi_k\psi_k\tbl\tbu \;, \label{T_3}
\end{eqnarray}
where $|\delta| < \epsilon$.
By using (\ref{T_3}) for the first three terms of (\ref{11_terms}), we get an
equation of the form (\ref{F_k}), where $F_k$ is given by the last six terms
in (\ref{11_terms}) plus the last six terms in (\ref{T_3}).
Terms involving the $\calD{\cdot}$ may be expressed in terms of $\tbl\tbu$
or $\tL\tU$
by using the definition of $V$, which in the current notation is
$$
v_{ij} = \frac{|\phi_i||\psi_j|}{\phi_i \psi_j}  \;.
$$
Consider the factor $\Phi\calD{6}(\psi_k\tbu)/\tr_{kk}$ in the term 
$-2\Phi\calD{6}(\psi_k\tbu)/\tr_{kk}$. We have
$$
(\Phi\calD{6}(\psi_k\tbu)/\tr_{kk})_{ij} = 
                      \phi_i\partial_j^{(6)}(\psi_k\tu_j)/\tr_{kk}\;.
$$
Recall that $\phi_i = [\phi_{i1}, \phi_{i2}]$ and $\psi_j = [\psi_{1j}, \psi_{2j}]$.
Then
$$
(\Phi\calD{6}(\psi_k\tbu)/\tr_{kk})_{ij} = 
    (\phi_{i1}\delta_{1j}^{(6)}\psi_{1k} + 
	\phi_{i2}\delta_{2j}^{(6)}\psi_{2k})\tu_j/\tr_{kk}\;,
$$
where $\delta_{1j}^{(6)}$ and $\delta_{1j}^{(6)}$ are the scaling factors from the
operator $\partial_j^{(6)}$. From the definition of $V$, using the fact that
$\tl_i \doteq \tr_{ik}/\tr_{kk}$, this can be shown to be
$$
(\Phi\calD{6}(\psi_k\tbu)/\tr_{kk})_{ij} =
   \hat{\delta}_{ij}^{(6)}v_{ik}b_{ik}^{-1}\tl_i\tu_j  \;,
$$
where $|\hat{\delta}_{ij}^{(6)}| \leq \max_{j=1,2}|\delta_{kj}^{(6)}|$\,.
In matrix form, we obtain
\beqn
\Phi\calD{6}(\phi_k\tbu)/\tr_{kk} = \hat{\Delta} \circ 
                                \diag\{v_{ik}/b_{ik}\} \tbl\tbu\;,
\eeqn
where $\hat{\Delta}$ and subsequent $\hat{\Delta}^{(\cdot)}$ 
are matrices with elements bounded in magnitude by $\epsilon$.
Similarly, all the other terms can be expressed as either 
\begin{itemize}
\item[(i)] an elementwise product
of $\hDel{\cdot}$ and a normal product of $\tbl\tbu$ 
and matrices derived from $B$ or $V$, or 
\item[(ii)] an elementwise product of the form 
$\hDel{\cdot} \circ B^I \circ V^{\prime} \circ L_{:,k+1:n}U_{k+1:n,:}$.
\end{itemize}
When this is done, the result follows.
\endproof

The next theorem uses Lemma~\ref{lem2} to obtain an elementwise bound for $|G|$.
\begin{theorem}
Let $H_k$ be as in Theorem~\ref{Th3.1}. Then
\begin{eqnarray}
|G| & \leq & c_1 b_{\min}^{-1}\hDel{1} \circ |\hat{L}||U|
   + c_2 b_{\min}^{-1}|L||\hat{U}| \circ \hDel{2}
   + c_3 b_{\min}^{-1}\hDel{3} \circ |L|\diag\{v_{kk}^{(k)}\}|U| \nonumber \\
 & + & c_4|B^I| \circ \hDel{4} \circ \sum_{k=2}^n |\hat{R}_k^{\prime}|
         \nonumber %
\end{eqnarray}
where $b_{\min}$ is the minimum modulus of the
elements of $B$, $\hat{L} = [\bv_{:k}^{(k)}]_{k=1}^n \circ L$,
$\hat{U} = U \circ [\bv_{k:}^{(k)}]_{k=1}^n$, and
$\hat{R}_k^{\prime} = V^{(k)} \circ L_{:,k:n}U_{k:n,:}$.
\end{theorem}
\proof
$G$ is evaluated by carrying out the summation in (\ref{lem2.1}), and using
the identities $\sum_{i=k}^n \ba_{:k}\bb_{k:} = AB$
and $\sum_{i=k}^n x_k\ba_{:k}\bb_{k:} = A\,\diag\{x_k\}B$.
\endproof

We now apply the last step in the above methodology to derive an
expression for $\|E\|$.

\pagebreak[4]
\begin{theorem}
\label{th:Ebound}
Let E be the backward error $E = \tL\tU - R$
in the factorization of $R$ using the GKO algorithm, let $\hat{L}$,
$\hat{U}$, $\hat R$, $B$ and $V$ be as above. Then $\|E\|$ is bounded by
\beq
\|E\| \leq \epsilon \left(c_5\frac{\bmax}{\bmin}g_1 + c_6ng_2\right) \|L\|\|U\| \;,
                                     \label{eq:Ebound}
\eeq
where the Frobenius norm is used,
$\bmax$ and $\bmin$ are the maximum and minimum moduli of the
elements of $B$, $c_5$ and $c_6$ are small constants, and $g_1$ and
$g_2$ are {\em generator growth factors}, defined by
\begin{eqnarray}
g_1 & = & c_7 \frac{\|\hat{L}\|}{\|L\|}
        + c_8 \frac{\|\hat{U}\|}{\|U\|}
        + c_9 \|\diag\{v_{kk}^{(k)}\}\| \;,             \Label{eq:g1}\\
g_2 & = & \max_{k=2,\ldots,n}\{ |\hat{R}_k\|/\|R_k\| \} \;,
                                                        \Label{eq:g2}
\end{eqnarray}
with $c_7, c_8, c_9 < 1$.
\end{theorem}
\proof From step 5 of the above methodology, we essentially invert the
Sylvester equation (\ref{dispE}) to derive an expression for $E$. To
do this we apply (\ref{coroll1.2}) in Corollary~\ref{coroll1}. This can
be written in matrix form
$$
E = B \circ G
$$
so
\begin{eqnarray}
|E| = |B| \circ |G|
    & \leq  & c_1 \bratio\hDel{5} \circ |\hat{L}||U|
            + c_2 \bratio\|L||\hat{U}| \circ \hDel{6} +   \nonumber  \\
    &       &
    	      c_3 \bratio\hDel{7} \circ |L|\diag\{v_{kk}\}|U| +
    	      c_4 \Del{8} \circ \sum_{k=2}^n |\hat{R}_k^{\prime}| \;.
           							\label{|E|}
\end{eqnarray}
We now define $g_2 \equiv \max_{k=2,\ldots,n}\|\hat{R}^{(k)}\|/\|R^{(k)}\|$, 
$g_4 \equiv \|\hat{L}\|/\|L\|$, $g_5 \equiv \|\hat{U}\|/\|U\|$ and\\ 
$g_6 \equiv \|\diag\{v_{kk}^{(k)}\}\|$.
These can be considered
to be generator growth factors --- they are functions of the $V^{(k)}$,
which from the definition (\ref{def_V}) are the ratio
of the products of the magnitudes of the generators to the products of the
generators. We will see in \S\ref{sec:Discussion} that these growth factors can
sometimes be large.

Taking the Frobenius norm of (\ref{|E|}), we can easily show that
\beq
\|E\| \leq  c_1 \delta_{1} \bratio g_3 \|L\|\|U\|
        + c_2 \delta_{2} g_4 \bratio\ \|L\|\|U\| 
        + c_3 \delta_{3}\bratio g_5 \|L\|\|U\| +    
          c_6 n \delta_{4} g_2 \|L\|\|U\| \;.    \label{||E||}
\eeq
where $0 \leq |\delta_{1}|,\ldots,|\delta_{4}| < \epsilon$.
The result follows by collecting the first three terms
of~(\ref{||E||}).~\endproof
The following corollary specializes the above result to the case when
$R$ is derived from a Toeplitz matrix.
\begin{corollary}
\label{coroll3.3}
Let $R$ be derived from a Toeplitz matrix  $T$ by the transformation
{\rm (\ref{T.C})} in\\ 
Theorem~\ref{Tgen.Cgen}, and let $c_1$, $c_2$, $g_1$, $g_2$ and $E$ be
as defined in Theorem~\ref{th:Ebound}. Then $\|E\|$ is bounded by
\beq
\|E\| \leq \epsilon c_{10} g_3 n \|L\|\|U\|   \label{||E(CToep)||}\;,
\eeq
where $c_{10} = \max(2c_5/\pi, c_6)$ and $g_3 = \max(g_1, g_2)$.
\end{corollary}
\proof Recall that $B = [1/(t_i-s_j)]$ is the ordinary Cauchy matrix
with displacement operator $\nabla_{\{D_s,D_t\}}$; 
from equations~(\ref{eq:Df})
in Theorem~\ref{Tgen.Cgen}, the $t_i$ are $n$ equally-spaced points
around the unit circle, including one at (1,0), and the $s_j$ are also
$n$ equally-spaced points around the unit circle, with each $s_j$
between two of the $t_i$. Clearly $\pi/n < t_i-s_j < 2 \;\; \forall i,j$,
so by the definition of $B$,
\beq
\bratio < 2n/\pi \;. \label{rhobound}
\eeq
Substituting (\ref{rhobound}) in (\ref{||E||}),
bounding $2c_5/\pi$ and $c_6$ by
$c_{10}$, and bounding $g_1$ and $g_2$ by $g_3$ yields the result.~\endproof

The above results show that the expressions for the backward error bounds
from the GKO algorithm are similar to the ones for Gaussian elimination
with partial pivoting (GE/PP)~\cite{G&VL},
except for the generator growth factors which
might arise in particular cases where the $\Phi^{(k)}$ and $\Psi^{(k)}$
are large, but not the $\Phi^{(k)}\Psi^{(k)}$ or the $R_k$. So there
may be some cases where large error growth may occur in the GKO
algorithm but not GE/PP. In \S\ref{sec:Discussion},
we give an example where this occurs.

\Section{Error Analysis of the GKO-Toeplitz Algorithm}
\label{sec:GKOTerror}

Recall that the steps in the GKO-Toeplitz algorithm are (i) compute the
generators from the Toeplitz matrix $T$ using (\ref{Phi.Toep}) and 
(\ref{Psi.Toep}), (ii) convert
them to generators of a Cauchy matrix using (\ref{T->C}) and 
(iii) compute factors 
$L$ and $U$ of this Cauchy matrix using the GKO algorithm.
The factors of $T$ are then given by (\ref{T_fac}). There are errors
incurred at each of these steps. In this section, we do not consider
permutations, as these do not contribute to the error. We will derive
a bound for the perturbation matrix $E_T$, defined by
\beqn
F^{*}\tL\tU FD = T + E_T \;.
\eeqn
In our development, we show in Theorem~\ref{th:ET_compts}
that $E_T$ consists of two components ---
the first due to the error $\|E\|$ incurred in the Cauchy factorization
and the second due to the errors incurred in computing the
Cauchy generators $\tilde{\Phi}$ and $\tilde{\Psi}$.
The latter is a Toeplitz-type perturbation $\Delta T$ such that
$T+\Delta T$ transforms exactly to $\tilde{\Phi}$ and $\tilde{\Psi}$. 
We then derive two lemmas needed to derive $\Delta T$, and then
present the main result of this section in Theorem~\ref{th:ETbound}.

\pagebreak[3]
\subsection{Main components of $E_T$}
$E_T$ has two main components, as is shown in the following.
\begin{theorem}
\label{th:ET_compts}
Let $F$ and $D$ be as in Theorem~\ref{Tgen.Cgen}, let $\tilde{\Phi}$
and $\tilde{\Psi}$ be the Cauchy generators computed using 
{\rm (\ref{Phi.Toep}), (\ref{Psi.Toep})} and {\rm (\ref{T->C})}, and 
let $\tL$ and $\tU$ be the factors computed from
$\tilde{\Phi}$ and $\tilde{\Psi}$ using the GKO algorithm.
Then the perturbed factorization of $T$ satisfies
\beq
F^{*}\tL\tU FD \equiv T + E_T = T - F^*E FD + \Delta T \;, \label{eq:ET_compts}
\eeq
where $E$ is as in Theorem~\ref{th:Ebound} and
$\Delta T$ is a Toeplitz-type perturbation of $T$ such that $T+\Delta T$ has
generators $\tilde{\Omega}$ and $\tilde{\Gamma}$ that transform exactly to
$\tilde{\Phi}$ and $\tilde{\Psi}$ using 
{\rm (\ref{Phi.Toep}),
(\ref{Psi.Toep})} and {\rm (\ref{T->C})}.
\end{theorem} 
\proof Let $\tR$ be the Cauchy matrix generated by $\tilde{\Phi}$ and 
$\tilde{\Psi}$. We have
$$
\tR = \tL\tU + E \;,
$$ 
and we know from (\ref{T.C}) that $\tilde{\Phi}$ and $\tilde{\Psi}$
are the generators for
$$
\tR = F(T+\Delta T)D^{-1}F^{*}
$$
where $T+\Delta T$ is some Toeplitz-type matrix. From the above two equations
we obtain
$$
T+\Delta T = F^{*}\tR FD = F^{*}(\tL\tU + E)FD \;,
$$
from which the desired result follows.
\endproof
 
\pagebreak[3]
Thus, by (\ref{eq:ET_compts}), we see that $E_T$ has one component with the
same norm bound as $E$, and another which perturbs $T$ to a matrix such that
its generators, say $\tilde{\Omega}$ and $\tilde{\Gamma}$, transform exactly 
to $\tilde{\Phi}$ and $\tilde{\Psi}$.
Before we derive an expression for $\Delta T$, we need two preliminary
results : expressions for $\tilde{\Omega}$ and $\tilde{\Gamma}$, and
a method to recover $T+\Delta T$ from its generators
$\tilde{\Omega}$ and $\tilde{\Gamma}$.

\subsection{Estimation of $\Delta T$ --- preliminary results}
The required results are given in the following two lemmas.
\begin{lemma}
\label{lem4.1}
Let $\Omega$ and $\Gamma$ be as in {\rm (\ref{Phi.Toep})} and 
{\rm (\ref{Psi.Toep})},
and let $\tilde{\Omega}$ and $\tilde{\Gamma}$ transform exactly 
to $\tilde{\Phi}$ and $\tilde{\Psi}$ using {\rm (\ref{T->C})}. 
Let $[\ba,\bb] = \tilde{\Omega} - \Omega$ and
let $[\bc,\bd\,] = \tilde{\Gamma}^* - \Gamma^*$. Then
\beqn
\ba = \bo
\eeqn
and $\|\bb\|$, $\|\bc\|$ and $\|\bd\|$ are bounded by
\begin{eqnarray}
\|\bb\| & \leq & \epsilon k_1 n^{3/2} \|\bomgii\| \;, \Label{bbb}\\
\|\bc\| & \leq & \epsilon k_2 n^{3/2} \|\bgami\| \;,  \Label{bbc}\\
\|\bd\| & \leq & \epsilon  \;.                        \Label{bbd}
\end{eqnarray}
\end{lemma}
\proof
We first consider the errors incurred in the computation of $\tphi$ and
$\tpsi$. We have
\begin{eqnarray*}
\tphi & = & fl\{\tF[\be_1, \tomgii]\}, 
       \quad\mbox{where}\;\; \tF = fl(F) \;, \tomgii = fl(\bomgii) \\
      & = & [\bi, fl(\tF \tomgii)] \;, \quad\mbox{where}\;\; 
                                    \bi = [1,1,\ldots,1]^T           \\
      & = & [\bi, \tF\tomgii + k_3 n \|\tomgii\|\del{1}]
\end{eqnarray*}
where $|\delta_i^{(1)}| < \epsilon\,, 
                 i=1,\ldots,n$.
After a few more steps, this becomes
\beqn
\tphi = F[\be_1, \bomgii+\bb]
\eeqn
where $\bb = \Del{7}\bomgii + k_4 (n+1) \|\tomgii\| F^{*}\del{1}$.
In a similar way, it can be shown that
\beqn
\tpsi^{*} = FD[\bgami^{*}+\bc, \be_n+\bd]
\eeqn
where $\bc = k_5 D^{*}\Del{8}D\bgami^{*} + k_6 (n+1) \|\bgami\|\del{2}$
and $\bd = D^{*}F^{*}d_n\Del{9}\boldf_{n:}^T$. 
Now the expressions in
square brackets transform exactly to $\tilde{\Omega}$ and $\tilde{\Gamma}$
respectively, and by taking norms of $\bb$, $\bc$ and $\bd$ the bounds
(\ref{bbb}) to (\ref{bbd}) can be demonstrated in a few steps.
\endproof
\begin{lemma}
\label{lem4.2}
For any matrix $A$, let $\nabla_{\{Z_1,Z_{-1}\}}A = B$. Then $A$ can be
recovered from $B$ using
\beq
a_{ij} = \sum_{k=j}^{n} b_{1+(i+k-j)\bmod\,n,k} 
         - \sum_{k=1}^{j-1} b_{1+(i+k-j)\bmod\,n,k}  \;.  \label{regen}
\eeq
\end{lemma}
\proof
{From} the displacement operator $\nabla_{\{Z_1,Z_{-1}\}}$, the following
properties of $B$ are easily seen:
\begin{eqnarray}
b_{ij} &=& a_{i-1,j} - a_{i,j-1}, \quad 1<i\leq n \;, 1\leq j<n \;,\label{B21} \\
b_{1j} &=& a_{nj} - a_{i,j+1}, \quad 1\leq j<n \;,                 \label{b11} \\
b_{in} &=& a_{i-1,j} + a_{i1}, \quad 1<i\leq n \quad\mbox{and}     \label{b22} \\
b_{1n} &=& a_{n,n-1} + a_{11} \;.                                  \label{b12}
\end{eqnarray}
It can be verified that if the elements of $A$ are given by
(\ref{regen}), then (\ref{B21}) to (\ref{b12}) are satisfied.
\endproof
Equation~(\ref{regen}) shows that an element $a_{ij}$ is recovered by computing
$x-y$, where $x$ is the sum of elements of $B$ down the diagonal, commencing
from $b_{i+1,j}$ and proceeding to the last column, wrapping from
the last row to the first if necessary during the summing;
$y$ is a similar ``wrapped diagonal sum'' from the first column to $b_{i,j-1}$.

\subsection{Main result}
We now use Theorem~\ref{th:ET_compts}, Lemma~\ref{lem4.1} and Lemma~\ref{lem4.2}
to derive a bound for the backward error $\|E_T\|$ in the GKO-Toeplitz
algorithm.
\begin{theorem}
\label{th:ETbound}
Let $F$ and $D$ be as in Theorem~\ref{Tgen.Cgen}, and
let $\tL$ and $\tU$ be the factors computed from $T$ using the
GKO-Toeplitz algorithm.
Then the perturbed factorization of $T$ satisfies
\beq
F^{*}\tL\tU FD \equiv T + E_T = T + E^{(1)} + E^{(2)}\;,	\label{E1E2}
\eeq
where $E^{(1)}$ is a general matrix with norm $\|E^{(1)}\| = \|E\|$,
$E$ is as in Theorem~\ref{th:Ebound}, and
$E^{(2)}$ is a Toeplitz-type matrix with norm bounded by
\beq
\|E^{(2)}\| \leq \epsilon c_{11} n^2 (\|\bt_{1:}\| + \|\bt_{:1}\|) \;.  \label{ETbound}
\eeq
\end{theorem}
\proof
By comparing (\ref{E1E2}) and (\ref{eq:ET_compts}), we see that $E^{(1)} = - F^*E FD$,
and because $F$ and $D$ are orthogonal matrices,
\beq
\|E^{(1)}\| = \|E\|  \;.                         \label{E1}
\eeq
{From} the above comparison we also have $E^{(2)} = \Delta T$,
a Toeplitz-type perturbation of $T$ such that $T+\Delta T$ has
generators $\tilde{\Omega}$ and $\tilde{\Gamma}$ that transform exactly to
the Cauchy generators
$\tilde{\Phi}$ and $\tilde{\Psi}$ computed using (\ref{Phi.Toep}). In the
following, we use $E^{(2)}$ for $\Delta T$. From Lemma~\ref{lem4.1}, we have
$$
\nabla(T+E^{(2)}) = \tilde{\Omega}\tilde{\Gamma}
              = \Omega\Gamma + \be_1\bc^* + \bomgii\bd^* + \bb\be_n^T \;,
$$
where $\bb$, $\bc$ and $\bd$ are bounded as in (\ref{bbb}) to (\ref{bbd}).
The second-order error term $\bb\bd^*$ has been omitted. We then have
$$
\nabla E^{(2)} = \be_1\bc^* + \bomgii\bd^* + \bb\be_n^T \;,
$$
and we use (\ref{regen}) to compute $E^{(2)}$. This yields, after some
algebra
$$
|\be_{:j}^{(2)}| = C_{j-1}(|\bc^R| + |\bb^R|) + |\bp_{:j}|
$$
where $C_k$ is a matrix which by premultiplication, circularly upshifts
a vector $k$ places, $\bx^R$ indicates the reversal of $\bx$, and the
moduli of $\bp_{:j}$ are bounded by
\begin{eqnarray}
|p_{ij}| &\leq& |\bomgii|^T C_{j-i-1} |\bd| \nonumber \\
         &\leq& \|\bomgii\| \|\bd\| \;.         \label{pbound}
\end{eqnarray}
Using (\ref{bbb}), (\ref{bbc}), (\ref{pbound}) and (\ref{bbd}) it is easily
seen that
$$
\|\be_{:j}^{(2)}\| \leq c_{12} n^{3/2} (\|\bomgii\| + \|\bgami\|)\;.
$$
{From} this, using the definitions (\ref{Phi.Toep}) and (\ref{Psi.Toep}),
we obtain the bound (\ref{ETbound}) for $E^{(2)}$. 
Together with (\ref{E1}), this yields the result.
\endproof

\Section{Discussion of Error Bounds}
\label{sec:Discussion}

We first discuss the factors in the above error bounds and relate them
to what would be expected for Gaussian elimination with
partial pivoting (GE/PP).
Then we show, for both the Cauchy and Toeplitz variants, that there are
some cases where the backward error growth can be large.

\subsection{Relation of bounds to those for GE/PP}
Consider the backward error $E$ incurred by the Cauchy variant
(equation~(\ref{eq:Ebound})). The term $\|L\|\|U\|$ is similar to that obtained
for GE/PP~\cite{G&VL}. However, the first factor contains
the generator growth factors $g_1$ and $g_2$. These are given by
ratios of norms of the hatted quantities to the unhatted quantities in
(\ref{eq:g1}) and (\ref{eq:g2}). The former are derived from the latter
by elementwise multiplication by submatrices of the $V^{(k)}$,  
which from their definitions (\ref{def_V}) are the ratio
of the products of the magnitudes of the generators to the products of the
generators. For an ordinary Cauchy matrix, $v_{ij}^{(k)} = 1 \; \forall i,j,k$
because $\Phi^{(k)}$ and $\Psi^{(k)}$ have only one column and row respectively.
However, for higher displacement-rank Cauchy matrices, there may be
significant cancellation in the computation of the denominator of (\ref{def_V}),
so they may be significant growth in the size of the $\hat{L}$, $\hat{U}$
and $\hat{R}_k$ compared to the $L$, $U$ and $R_k$ respectively.

The backward error $E_T$ incurred by the Toeplitz variant has two
components --- one with the same norm as $E$ above, and a Toeplitz-type
component with norm bounded as in (\ref{ETbound}). The latter bound
is proportional to $n^2$ and contains no growth factors, so it would
be expected that the bound would be dominated by the first component.

We next give examples where the generator growth might be expected to
be large in the Cauchy and Toeplitz variants.

\subsection{Examples of large generator growth}

{\em Cauchy case}.
\ Here, we can select an example where all the elements of $V = V^{(1)}$
are large. This will occur when significant cancellation occurs in
the computation of the $\phi_i\psi_j$. Such an example is
$$
\Phi = [\ba, \ba+\boldf] \;, \quad \Psi = [\ba, -\ba]^T \;, %
$$
where $\|\ba\|$ is of order unity, and $\|\boldf\|$ is very small. Then
$\Phi\Psi = -\boldf\ba$, that is, all the elements of $\Phi\Psi$ are
very small compared to those of $|\Phi||\Psi|$. Moreover, because
$\ba$ and $\boldf$ can be arbitrary except for their norms, the
original matrix $[(t_i-s_j)^{-1} \phi_i\psi_j]$ is in general well-conditioned.

{\em Toeplitz case}.
The Toeplitz case has an extra constraint on the selection of $\Phi$ and
$\Psi$, since it must be generated from $\Omega$ and $\Gamma$ using
the transformations (\ref{T->C}). Because of this constraint, there
is no case where all the elements of $V$ can be made large. However,
all of the first column of $V$ can be made large, and this will cause
error growth, in spite of the pivoting. This will happen in the
following case.

Recall that $a_{i-j} = t_{ij}\; \forall i,j$. Select 
\begin{eqnarray}
                 a_0 &=& 1                              \label{a_0}  \\
\mbox{and} \quad a_i &=& -a_{i-n} \;,\quad 1\leq i\leq n-1 \;, \label{a_i}
\end{eqnarray} 
so that $\Omega = [\be_1, \be_1]$.
Then all of the first column of $V$ will be large if
$\psi_{11}+\psi_{12}$ is very small compared to $\psi_{11}$
and $\psi_{12}$. It can be verified from (\ref{T->C}) that if we select 
$a_1,\ldots,a_{n-1}$ to satisfy
\beq
\sum_{j=1}^{n-1} a_{n-j}\exp(i\pi(j-1)/n) = -\exp(i\pi(n-1)/n + \delta/2
                   \label{delta_cond}
\eeq
then $\psi_{11}+\psi_{12} = \delta$. There is a wide variety of choices
for the $a_j$. Let $n$ be even, and set 
\beq
a_{n-j}=0                \label{anmj}
\eeq
except for $a_{n/2-1}$ and $a_{n-1}$. Then (\ref{delta_cond}) is satisfied
when 
\beq
a_{n/2-1} = -\sin(\pi/n) + \Im(\delta/2) \;, \quad
a_{n-1} = \cos(\pi/n) + \Re(\delta/2)  \;.     \label{anm1}
\eeq
So if $\delta$ is small, and the $a_j$ are selected according to
(\ref{a_0}), (\ref{anmj}), (\ref{anm1}) and (\ref{a_i}), all of the
first column of $V$ will be large, with magnitude $O(1/\delta)$.

{\em Numerical examples}.
\ Order-8 Toeplitz matrices were generated according to
(\ref{a_0}), (\ref{anmj}), (\ref{anm1}) and (\ref{a_i}), with
\mbox{$\delta = 10^{-k}$}, $k=2,\ldots,16$. For each matrix, the system
$T\bx = \bi$ was solved. It was found that the normalized solution
error $\|\tilde{\bx}-\bx\|/\|\bx\|$ grew as the square of $1/\delta$,
and the normalized residual $\|T\tilde{\bx}-\bi\|/\|\bb\|$ grew
linearly with $1/\delta$. Thus the algorithm is only weakly stable
in this case.

\Section{Modified GKO Algorithm}
\label{sec:GKOmod}

The problem with the original pivoting strategy is that when all elements of
$\br_{:1}$ are small and all elements of $\bv_{:1}$ are large, normal
partial pivoting will not stabilize the algorithm.
Complete pivoting will
do so, but requires $O(n^2)$ operations to find the pivot at each major
step and $O(n^3)$ operations overall. However, a strategy of using the
largest element in the first row {\em and column} should stabilize the
algorithm in most cases, and we see that it does in the above cases.

\pagebreak[3]
To incorporate this row-1/column-1 pivoting, it is easy to see
that the following steps should be added to the GKO algorithm
(Algorithm~\algGKOC):
\begin{itemize}
   \item Step 1: add substep $P^{\prime} \leftarrow I$, where $P^{\prime}$
         will be the matrix of column interchanges.
   \item After loop to recover column 1 of $R_k$ : add loop to recover
         row 1 of $R_k$.
   \item After loop to find maximum $\max_1$ in column~1 : add loop to find 
         maximum $\max_2$ in row~1. If $\max_1 \geq \max_2$, carry out
         row interchanges as in Algorithm~\algGKOC. Otherwise carry out
         column interchanges by swapping the appropriate elements
         in $\bs$, $\Psi^{(k)}$ and $\br_{k:}^{(k)}$, and the appropriate
         columns in $U$ and $P^{\prime}$.
   \item After computation of $L$, $U$, $P$ and $P^{\prime}$, the
         factors of $R$ are $P^TLUP^{\prime T}$.
\end{itemize}

{\em Results}.
\ When the modified algorithm was used on the same set of systems as was
considered in the previous section, it was found that the normalized solution
error $\|\tilde{\bx}-\bx\|/\|\bx\|$ grew linearly with $1/\delta$
and the condition number of $T$,
and the normalized residual $\|T\tilde{\bx}-\bi\|/\|\bb\|$ was approximately
constant at about $4 \times 10^{-15}$, a small multiple of $\epsilon$.
Thus the modified algorithm is stable in this case.

\Section{Conclusions}
\label{sec:Conclusion}

It has been shown that bound for the backward error in the GKO algorithm
is similar to that for partial pivoting, except that extra factors,
the generator growth factors, are included.
These factors can be large when there
is sufficient cancellation in the computation of the generators. Examples
of this have been presented, and it was demonstrated that the original
GKO algorithm was only weakly stable in these cases. A modified version
which uses row~1/column~1 pivoting was then presented; this version was
stable in these cases.

It is not known whether there are any cases upon which the modified
algorithm will give large errors. Further work needs to be done to
ascertain this, and if such cases can be found, the pivot strategy needs
to be improved further. The aim is to find the maximum in $R$, or an
element close to the maximum, still in $O(n)$ operations. An extension
of the above strategy may be to have a few iterations in the search,
i.e. search for the row-1/column-1 maximum, say at $r_{1p}$, then
search along column $p$ for the maximum there, and so on. This may find
a better pivot at the expense of some extra work.

A practical strategy is to use the modified algorithm of \S\ref{sec:GKOmod}
followed by a check of the residual; in the unlikely event that the residual
is large we can resort to a stable $O(n^3)$ algorithm.

\end{document}